\documentclass[12pt,twoside]{article}
\usepackage[height=22cm,width=14cm]{geometry}
\usepackage{amssymb,amsmath,amsthm,latexsym}
\usepackage{amsfonts}
\usepackage{graphicx}
\usepackage{subcaption}
\usepackage{mathtools}
\usepackage{longtable}
\usepackage{makecell}
\usepackage{cite}
\usepackage{lipsum}
\usepackage{xcolor}
\usepackage{ragged2e}

\newtheorem{theorem}{Theorem}[section]

\vspace{5cm}

\begin{document}
	
	\label{'ubf'}
\setcounter{page}{1}

\begin{center}
{ 
       {\Large \textbf { On the Elliptic Sombor and Euler Sombor indices of Corona product of certain graphs
                               }
       }
\\

\medskip

{\sc B. Kirana$^{1}$, M.C. Shanmukha$^{2,*}$,  A. Usha$^{3}$ }\\

{\footnotesize $^{1}$ Department of Mathematics, KVG College of Engineering, Sullia, 574327  and affiliated to Visvesvaraya  Technological University, Belagavi-590018, India}\\

{\footnotesize $^{2}$ Department of Mathematics, PES Institute of Technology and Management, Shivamogga, 577204 and  affiliated to Visvesveraya  Technological University, Belagavi, 590018, India}\\

{\footnotesize $^3$Department of Mathematics, Alliance School of Applied Mathematics, Alliance University, Bangalore, 562106, India}\\

{\footnotesize e-mail*: {\it mcshanmukha@gmail.com}}}
\end{center}
\thispagestyle{empty}

\hrulefill\begin{abstract}
Elliptic Sombor and Euler Sombor indices are recently defined topological indices using Sombor index. Elliptic sombor index is defined as $ESO(G)=\sum_{uv\in E(G)}(d_u+ d_v)\sqrt{d^2_u+ d^2_v}$ and Euler Sombor index is defined as 
$EU(G)= \sum_{uv\in E(G)}\sqrt{{d_{u}^2+d_{v}^2}+d_ud_v}$, where $d_u$ and $d_v$ are degrees of vertices $u$ and $v$ in graph $G$. In this article, we compute the elliptic Sombor and Euler Sombor indices of some resultant graphs. Using the operations join and Corona product on standard graphs like path, cycle and complete graphs.
\end{abstract}

\hrulefill\\


\section{Introduction}

Degree-based topological indices are fundamental tools in the fields of mathematical chemistry and graph theory, serving as quantitative measures of a molecular graph structure. A molecular graph represents a chemical compound where vertices correspond to atoms and edges correspond to chemical bonds. The degree of a vertex in this context refers to the number of edges connected to it, representing the number of bonds of an atom. The primary utility of these indices lie in their ability to predict various physical, chemical, and biological properties of molecules\cite{H, NT}.

The Sombor index, introduced by Gutman in 2021, is a topological index derived from graph theory and Euclidean geometry. This index has gained significant traction in the fields such as chemistry and pharmacology for its utility in characterizing molecular structures and predicting biological activity. Let us delve into the formal definition, properties, and an example calculation of the Sombor index \cite{GI2, MI}. 

The general form of a vertex degree-based index is a function which is chosen such that it satisfies symmetry property. The edge $uv$ representation in 2-dimensional coordinate system is called the degree-point of the edge $uv$. The Euclidean distance between the degree-point $(d_u,d_v)$ and the origin $O$ is $\sqrt{d_u^2+d_v^2}$ which is the definition of Sombor index\cite{IG3}.

In 2023, Gutman et al., \cite{GI3} introduced another version of Sombor index called elliptic Sombor index referring to the orbits of planets in the solar system which takes elliptic orbits with the sun as focus point. The concept of an elliptic Sombor index can be an extension of the traditional Sombor index by incorporating elliptic distance rather than the Euclidean distance between degree-points.

In 2024, Gutman et al., \cite{IG} showed that the lengths of semi major and the semi minor axes are equal in an ellipse. Leonard Euler found the approximate perimeter of the ellipse as $\pi \sqrt{2(d_u^2+d_v^2)(d_u+d_v)^2}$. Using these relations, Euler Sombor index was proposed as $\sqrt{d_u^2+d_v^2+d_u.d_v}$. Algebraically, there is a geometric analogy of Sombor and Euler Sombor indices.

Bibhas et al., \cite{B} studied structural properties of corona graphs including the statistics of signed links and types of signed triangles and degree distribution. They analyzed algebraic conflict of signed corona graphs generated by seed graphs.  Sheeba et al., \cite{SH} established  the explicit expressions of Y -index of different types of corona product of graphs. Khalid et al., \cite{Kh} studied certain degree-based topological indices such as Randić index, Zagreb indices, multiplicative Zagreb indices, Narumi–Katayama index, atom-bond connectivity index, augmented Zagreb index, geometric-arithmetic index, harmonic index, and sum-connectivity index for the bistar graphs and the corona product. Sheeja et al., \cite{She} established the explicits expressions for the SK index over different types of corona products on graphs are presented. Iswadi et al., \cite{RG} determined the metric dimension of corona product, and lower bound of metric dimension of join of graphs. Dhananjaya et al., \cite{Kum} derived ABC index and the GA index of several corona products of graphs composed of path, cycle, and complete graphs.
Arif et al., \cite{AK1} computed the Sombor index of some graph operations namely, join and corona product of two graphs, for the standard graphs path, cycle, and complete graphs. 

 Motivated by the above studies on graph operations, this study focusses on computing the elliptic Sombor and Euler Sombor indices of some resultant graphs using the operations join and Corona product on standard graphs like path, cycle and complete graphs.

Graph operations, such as join, union, Cartesian product, and others, are fundamental tools in graph theory that significantly influence topological indices \cite{F, MV, GW, AK2}. Given two disjoint graphs \( G_1 = (V_1, E_1) \) and \( G_2 = (V_2, E_2) \), the join \( G_1 + G_2 \) (See Figure 1) is defined as:
\[
G_1 + G_2 = (V_1 \cup V_2, E_1 \cup E_2 \cup \{(v_1, v_2) \mid v_1 \in V_1, v_2 \in V_2\})
\]
Let $G$ be a connected graph with $n$ vertices and $H$ be a graph with at least two vertices, the corona product of $G$ and $H$, is defined as a graph which is formed by taking $n$ copies of $H$ and connecting $i^{th}$ vertex of $G$ to the vertices of $H$ in each copies, and it is denoted by $G\odot H$ (See Figure 2).
\begin{figure}[ht!]
    \centering
    \renewcommand\thefigure{1}
    \includegraphics[height=5cm, width=7cm]{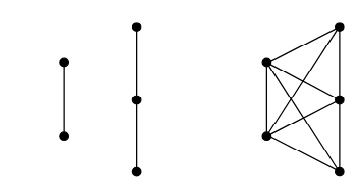}\\
    \label{fi:twentytwo}
    \caption{Join of two paths $P_2$ and $P_3$.}
    \end{figure}
    \begin{figure}[hbt!]
    \centering
    \renewcommand\thefigure{2}
    \includegraphics[height=5cm, width=7cm]{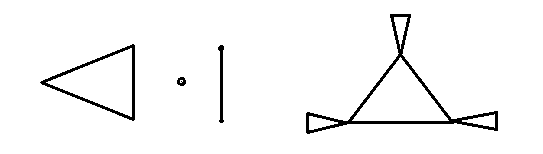}\\
    \label{fi:twentytwo}
    \caption{Corona product of two graphs $K_3$ and $P_2$.}
    \end{figure}\\
Elliptic Sombor and Euler Sombor indices are defined as follows\cite{IG,GI3}
\begin{align}
& ESO(G)=\sum_{uv\in E(G)}(d_u+ d_v)(\sqrt{d^2_u+ d^2_v})\\
& EU(G)= \sum_{uv\in E(G)}\sqrt{{d_{u}^2+d_{v}^2}+d_ud_v}
\end{align}
\section{ Main Results}

In this section, Elliptic Sombor and Euler Sombor indices of $P_r+P_s$, $C_r+C_s$, $K_r+K_s$, $P_r \odot P_s$, $C_r \odot C_s$ are obtained.

If $P_r$ and $P_s$ are two paths of order $|V(G)|=r+s$ and size $|E(G)|=(r-1)+(s-1)+rs=r+s+rs-2$.
In the following theorem, the Elliptic Sombor and Euler Sombor indices of two path graphs $P_r$ and $P_s$ are computed.
\begin{theorem}
Elliptic Sombor index of $P_r+P_s$ is given by 
\begin{alignat*}{2}
ESO(P_r+P_s)=\begin{cases}
 108\sqrt{2},\, r=s=2\\\\
70+2\sqrt{2}(s+1)^2+32\sqrt{2}(s-3)\\
+4(s+4)\sqrt{9+(s+1)^2}\\
+2(s-2)(s+5)\sqrt{16+(s+1)^2}; s>2, r=2\\\\
70+2\sqrt{2}(r+1)^2+32\sqrt{2}(r-3)\\
+4(r+4)\sqrt{9+(r+1)^2}\\
+2(r-2)(s+5)\sqrt{16+(r+1)^2};\, r>2, \,s=2\\\\
2(2s+3)\sqrt{(s+1)^2+(s+2)^2}\\
+(r-3)(2s+4(s+2)\sqrt{2}+2(2r+3)\sqrt{(r+1)^2+(r+2)^2}\\
+(s-3)(2r+4)(r+2)\sqrt{2}+4(r+s+2)\sqrt{(r+1)^2+(s+1)^2} ; r,s>2
\end{cases}
\end{alignat*}
\begin{proof}
Consider the join of two graphs $P_r$ and $P_s$. If $r,s>1$,  based on the degrees of vertices of
$P_r+P_s$, there are 4 types of vertices having degrees $s+1, s+2, r+1$and $r+2$. Then by
considering the degrees of the vertices, there are 10 types of edge partitions of $P_r+P_s$ as shown in the Table 1.
\begin{table}[h!]
\centering
\caption{Number of edges for different cases in $P_r + P_s$.}
\begin{tabular}{|c|c|c|}
\hline
\textbf{Case} & \textbf{Number of edges} & \textbf{Condition} \\
\hline
$(s + 1, s + 1)$ & $1$ & $r = 2$ \\
                 & $0$ & $r > 2$ \\
\hline
$(s + 1, s + 2)$ & $0$ & $r = 2$ \\
                 & $2$ & $r > 2$ \\
\hline
$(s + 2, s+ 2)$ & $0$ & $r = 2$ \\
                 & $r - 3$ & $r> 2$ \\
\hline
$(r+ 1, r+ 1)$ & $1$ & $s= 2$ \\
                 & $0$ & $s > 2$ \\
\hline
$(r + 1, r + 2)$ & $0$ & $s = 2$ \\
                 & $2$ & $s > 2$ \\
\hline
$(s + 2, s + 2)$ & $0$ & $s = 2$ \\
                 & $s - 3$ & $s> 2$ \\
\hline
$(s + 1, r+ 1)$ & $4$ & \text{Always} \\
\hline
$(s+ 1, r+ 2)$ & $0$ & $s = 2$ \\
                 & $2(s - 2)$ & $s> 2$ \\
\hline
$(s + 2, s + 1)$ & $0$ & $r= 2$ \\
                 & $2(r - 2)$ & $r> 2$ \\
\hline
$(s+ 2, r + 2)$ & $0$ & $r = 2$ or $s = 2$ \\
                 & $(r- 2)(s- 2)$ & $r, s > 2$ \\
\hline
\end{tabular}
\label{table:edges}
\end{table}\\\\
\textbf{Case 1}\,\, r=s=2
\begin{alignat*}{2}
ESO(P_r+P_s)=&2(s+1)\sqrt{(s+1)^2+(s+1)^2}+2(r+1)\sqrt{(r+1)^2+(r+1)^2},\\
&+4(s+1+r+1)\sqrt{(s+1)^2+(r+1)^2}\\
ESO(P_r+P_s)=&2(s+1)\sqrt{2(s+1)^2}+2(r+1)\sqrt{2(r+1)^2}\\
&+4(r+s+2)\sqrt{(s+1)^2+(r+1)^2},\\
\text{Substitute}\, r=2,s=2\, \text{we get}\\
ESO(P_r+P_s)=&108\sqrt{2}.
\end{alignat*}
\textbf{Case 2} r=2, s$>$2\\
\begin{alignat*}{2}
ESO(P_r+P_s)=&2(s+1)\sqrt{2(s+1)^2}+2(r+1+r+2)\sqrt{(r+1)^2+(r+2)^2}\\
&+(s-3)2(r+2)\sqrt{2(r+2)^2}+4(r+1+s+1)\sqrt{(r+1)^2+(s+1)^2}\\
&+2(s-2)(r+2+s+1)\sqrt{(r+2)^2+(s+1)^2},\\
ESO(P_r+P_s)=&2(s+1)^2\sqrt{2}+2(2r+3)\sqrt{(r+1)^2+(r+2)^2}\\
&+2(s-3)(r+2)(r+2)\sqrt{2}+4(r+s+2)\sqrt{(r+1)^2+(s+1)^2}\\
&+2(s-2)(r+s+3)\sqrt{(r+2)^2+(s+1)^2},
\end{alignat*}
take $r=2$ we obtain
\begin{alignat*}{2}
ESO(P_r+P_s)=&2\sqrt{2}(s+1)^2+2\times 7\sqrt{3^2+4^2}+32(s-3)\sqrt{2}\\
&+4(s+4)\sqrt{9+(s+1)^2}+2(s-2)(s+5)\sqrt{16+(s+1)^2},\\
ESO(P_r+P_s)=&70+2\sqrt{2}(s+1)^2+32(s-3)\sqrt{2}+4(s+4)\sqrt{9+(s+1)^2}\\
&+2(s-2)(s+5)\sqrt{16+(s+1)^2}.
\end{alignat*}
\textbf{Case 3}  s=2, r$>$2\\
\begin{alignat*}{2}
ESO(P_r+P_s)=&2(r+1)\sqrt{2(r+1)^2}+2(s+1+s+2)\sqrt{(s+1)^2+(s+2)^2}\\
&+(r-3)2(s+2)\sqrt{2(s+2)^2}+4(r+1+s+1)\sqrt{(r+1)^2+(s+1)^2}\\
&+2(r-2)(r+1+s+2)\sqrt{(r+1)^2+(s+2)^2},\\
ESO(P_r+P_s)=&2(r+1)^2\sqrt{2}+2(2s+3)\sqrt{(s+1)^2+(s+2)^2}\\
&+2(r-3)(s+2)^2\sqrt{2}+4(r+s+2)\sqrt{(r+1)^2+(s+1)^2}\\
&+2(r-2)(r+s+3)\sqrt{(s+2)^2+(r+1)^2},\\
\text{Substitute} \,\,s=2\\
ESO(P_r+P_s)=&2\sqrt{2}(r+1)^2+2\times 7\sqrt{3^2+4^2}+32(r-3)\sqrt{2}\\
&+4(r+4)\sqrt{9+(r+1)^2}+2(r-2)(s+5)\sqrt{16+(r+1)^2},\\
ESO(P_r+P_s)=&70+2\sqrt{2}(r+1)^2+32\sqrt{2}(r-3)+4(r+4)\sqrt{9+(r+1)^2}\\
&+2(r-2)(s+5)\sqrt{16+(r+1)^2}.
\end{alignat*}
\textbf{Case 4}\,\, r, s$>$2
\begin{alignat*}{2}
ESO(P_r+P_s)=&2(2s+3)\sqrt{(s+1)^2+(s+2)^2}+2(r-3)(s+2)\sqrt{2(s+2)^2}\\
&+2(2r+3)\sqrt{(r+1)^2+(r+2)^2}+2(s-3)(r+2)\sqrt{2(r+2)^2}\\
&+4(r+s+2)\sqrt{(r+1)^2+(s+1)^2}+2(s-2)(r+s+3)\\
&\sqrt{(r+2)^2+(s+1)^2}+(r-2)(s-2)(r+s+4)\sqrt{(r+2)^2+(s+2)^2}\\
&+2(r-2)(r+s+3)\sqrt{(r+1)^2+(s+2)^2}.
\end{alignat*}
\end{proof}
\end{theorem}
\begin{theorem}
The Euler Sombor index of $P_r+P_s$ is given by 
\begin{alignat*}{2}
EU(P_r+P_s)=\begin{cases}
 18\sqrt{3}, r=s=2\\\\ 
2\sqrt{37}+(5s-11)\sqrt{3}+4\sqrt{(s+1)^2+3s+12}\\
+2(s-2)\sqrt{(s+1)^2+4s+20}; s>2, r=2\\\\
2\sqrt{37}+(5s-11)\sqrt{3}+4\sqrt{(r+1)^2+3r+12}\\
+2(r-2)\sqrt{(r+1)^2+4r+20}; r>2, s=2\\\\
2\sqrt{(s+1)^2+(s+2)^2+(s+1)(s+2)}+(r-3)\sqrt{3(s+2)^2}\\
+2\sqrt{(r+1)^2+(r+2)^2+(r+1)(r+2)}+(s-3)\sqrt{3(r+2)^2}\\
+4\sqrt{(r+1)^2+(s+1)^2+(r+1)(s+1)}\\
+2(s-2)\sqrt{(r+2)^2+(s+1)^2+(r+2)(s+1)}\\
+(r-2)(s-2)\sqrt{(r+2)^2+(s+2)^2+(r+2)(s+2)}\\+2(r-2)\sqrt{(r+1)^2+(s+2)^2+(r+1)(s+2)}; r,s>2
\end{cases}
\end{alignat*}
\begin{proof}
 Using edge partition of join of two path graphs $P_r$ and $P_s$ from Table 1,
 Euler Sombor index is calculated using the equation 2. The following cases are considered.\\
 \textbf{Case 1} r=s=2\\
\begin{alignat*}{2}
EU(P_r+P_s)=&\sqrt{(s+1)^2+(s+1)^2+(s+1)(s+1)}\\
&+\sqrt{(r+1)^2+(r+1)^2+(r+1)(r+1)}\\
&+4\sqrt{(s+1)^2+(r+1)^2+(s+1)(r+1)},
\end{alignat*}
\begin{alignat*}{2}
EU(P_r+P_s)=&\sqrt{3(s+1)^2}+\sqrt{3(r+1)^2}+4\sqrt{(s+1)^2+(r+1)^2+(s+1)(r+1)},\\
take\, r=2,s=2\, \, we\, get\\
EU(P_r+P_s)=&18\sqrt{3}.
\end{alignat*}
\textbf{Case 2} r=2, s$>$2\\
\begin{alignat*}{2}
EU(P_r+P_s)=&\sqrt{2(s+1)^2+(s+1)(s+1)}+2\sqrt{(r+1)^2+(r+2)^2+(r+1)(r+2)}\\
&+(s-3)\sqrt{2(r+2)^2+(r+2)(r+2)}\\
&+4\sqrt{(r+1)^2+(s+1)^2+(r+1)(s+1)}\\
&+2(s-2)\sqrt{(r+2)^2+(s+1)^2+(r+2)(s+1)},\\
EU(P_r+P_s)=&\sqrt{3}[(s+1)+(r+2)(s-3)]+2\sqrt{(r+1)^2+(r+2)^2+(r+1)(r+2)}\\
&+4\sqrt{(r+1)^2+(s+1)^2+(r+1)(s+1)}\\
&+2(s-2)\sqrt{(r+2)^2+(s+1)^2+(r+2)(s+1)},\\
\text{Substitute} \, r=2\\ 
EU(P_r+P_s)=&2\sqrt{37}+(5s-11)\sqrt{3}+4\sqrt{(s+1)^2+3s+12}+2(s-2)\sqrt{(s+1)^2+4s+20}.
\end{alignat*}
\textbf{Case 3} $s=2$, $r>2$\\
\begin{alignat*}{2}
EU(P_r+P_s)=&\sqrt{2(r+1)^2+(r+1)(r+1)}+2\sqrt{(s+1)^2+(s+2)^2+(s+1)(s+2)}\\
&+(r-3)\sqrt{3(r+2)^2}+4\sqrt{(r+1)^2+(s+1)^2+(r+1)(s+1)}\\
&+2(r-2)\sqrt{(s+2)^2+(r+1)^2+(s+2)(r+1)},\\
EU(P_r+P_s)=&\sqrt{3}[(r+1)+(s+2)(r-3)]+2\sqrt{(s+1)^2+(s+2)^2+(s+1)(s+2)}\\
&+4\sqrt{(r+1)^2+(s+1)^2+(r+1)(s+1)}\\
&+2(r-2)\sqrt{(s+2)^2+(r+1)^2+(r+2)(s+1)},\\
\text{Substitute} \, s=2 \\
EU(P_r+P_s)=&2\sqrt{37}+(5s-11)\sqrt{3}+4\sqrt{(r+1)^2+3r+12}\\
&+2(r-2)\sqrt{(r+1)^2+4r+20}.
\end{alignat*}
\textbf{Case 4 } $r,s>2$
\begin{alignat*}{2}
 EU(P_r+P_s)=&2\sqrt{(s+1)^2+(s+2)^2+(s+1)(s+2)}+(r-3)\sqrt{3(s+2)^2}\\
&+2\sqrt{(r+1)^2+(r+2)^2+(r+1)(r+2)}+(s-3)\sqrt{3(r+2)^2}\\
&+4\sqrt{(r+1)^2+(s+1)^2+(r+1)(s+1)}\\
&+2(s-2)\sqrt{(r+2)^2+(s+1)^2+(r+2)(s+1)}\\
&+(r-2)(s-2)\sqrt{(r+2)^2+(s+2)^2+(r+2)(s+2)}\\
&+2(r-2)\sqrt{(r+1)^2+(s+2)^2+(r+1)(s+2)}.
\end{alignat*} 
\end{proof}
\end{theorem}
 \begin{theorem}
 Let $r,s$ be two positive integers such that $r,s\geq 3$. Then Elliptic Sombor index and Euler Sombor index of $C_r+ C_s$ are given by 
 \begin{alignat*}{2}
 &ESO(C_r+C_s)=2\sqrt{2}(s+2)^2+s(r+2)^2]+rs(r+s+4)\sqrt{(r+2)^2+(s+2)^2}.\\
 &EU(C_r+C_s)=2\sqrt{3}(r+s+rs)+rs\sqrt{(r+2)^2+(s+2)^2+(r+2)(s+2)}.
 \end{alignat*}
 \begin{proof}
 Consider join of two cycles $C_r$ and $C_s$ which is graph with order $|V(C_r+C_s)|=r+s$ and size $|E(C_r+C_s)|=r+s+rs$. There are two types of degree vertices such as $s+2$ and other type has $r+2$  in the resultant graph. Thus, there are three edge partitions
\begin{table}[h!]
\caption{Edge partitions of $C_r + C_s$.}
\centering
\begin{tabular}{|c|c|}
\hline
\textbf{Case} & \textbf{Number of edges}\\
\hline
$(s +2 , s +2)$ & $r$ \\
\hline
$(s+2, r + 2)$ & $rs$ \\
\hline
$(r + 2, r + 2)$ & $s$\\
\hline
\end{tabular}
\label{table:edges}
\end{table}\\\\
Using above edge partitions and the formula 
\begin{alignat*}{2}
ESO(C_r+C_s)=&r(s+2+s+2)\sqrt{2(s+2)^2}+rs(r+2+s+2)\sqrt{(r+2)^2+(s+2)^2}\\
&+s(r+2+r+2)\sqrt{2(r+2)^2},\\
ESO(C_r+C_s)=&2\sqrt{2}r(s+2)^2+rs(r+s+4\sqrt{(r+2)^2+(s+2)^2}+2\sqrt{2}s(r+2)^2,\\
ESO(C_r+C_s)=&2\sqrt{2}\left((r(s+2)^2+s(r+2)^2\right)+rs(r+s+4)\sqrt{(r+2)^2+(s+2)^2}.
\end{alignat*}
\begin{alignat*}{2}
EU(C_r+C_s)=&r\sqrt{3(s+2)^2}+rs\sqrt{(r+2)^2+(s+2)^2+(r+2)(s+2)}+s\sqrt{3(r+2)^2},\\
EU(C_r+C_s)=&\sqrt{3}\left[r(s+2)+s(r+2)\right]+rs\sqrt{(r+2)^2+(s+2)^2+(r+2)(s+2)},\\
EU(C_r+C_s)=&2\sqrt{3}\left[r+s+rs\right]+rs\sqrt{(r+2)^2+(s+2)^2+(r+2)(s+2)}.
\end{alignat*}
 \end{proof}
 \end{theorem}
 \begin{theorem}
 If two complete graphs $K_r$ and $K_s$ are joined, the resultant is a complete graph $K_{r+s}$ with order $|V(K_{r+s})|=r+s$ and size $|E(K_{r+s})|=\frac{(r+s)(r+s-1)}{2}$. Then \\
 \begin{alignat*}{2}
 ESO(K_r+K_s)=\sqrt{2}(r+s)(r+s-1)^3.\\
 EU(K_r+K_s)=\frac{\sqrt{3}}{2}(r+s)(r+s-1)^2.
 \end{alignat*}
 \begin{proof}
Each vertex of join of two graphs $K_r$ and $K_s$ has degree $r+s-1$. So there is only one type of edge that is $(r+s-1, r+s-1)$. Using equation (1), Elliptic Sombor index is 
\begin{alignat*}{2}
ESO(K_r+K_s)=&ESO(K_{r+s})\\
=&\frac{(r+s)(r+s-1)}{2}(r+s-1+r+s-1)\sqrt{2(r+s-1)^2}\\
=&\sqrt{2}{(r+s)(r+s-1)^3}.\\
\end{alignat*}
Using equation (2), Euler Sombor index is\\
\begin{alignat*}{2}
EU(K_r+K_s)=&EU(K_{r+s})\\
=&\frac{(r+s)(r+s-1)}{2}\sqrt{3(r+s-1)^2}\\
=&\sqrt{3}\frac{(r+s)(r+s-1)^2}{2}.
\end{alignat*}
\end{proof}
\end{theorem}
 \begin{theorem}
 Let $C_r$ and $K_S$ be the cycle and complete graph respectively. Then Elliptic Sombor index and Euler sombor index of $C_r+K_s$ are
 \begin{alignat*}{2}
ESO(C_r+K_s)=&2\sqrt{2}r(s+2)^2+rs(r+2s+1)\sqrt{(s+2)^2+(r+s-1)^2}\\
&+C(s,2)(r+s-1)\sqrt{2}.\\
EU(C_r+K_s)=&r(s+2)\sqrt{3}+rs\sqrt{(s+2)^2+(r+s-1)^2+(s+2)(r+s-1)}\\
&+C(s,2)(r+s-1)\sqrt{3}.
\end{alignat*}
\begin{proof}
The graph  $C_r+K_s$ contains $r+s$ number of vertices and $r+rs+C(s,2)$ number of edges. Here $C_r$ has $r$ number of edges and $K_s$ has $C(s,2)$ number of edges. There are three types of edge partitions that are tabulated in Table 3. 
\begin{table}[h!]
\caption{Number of edges for different cases in $C_r + K_s$.}
\centering
\begin{tabular}{|c|c|}
\hline
\textbf{Case} & \textbf{Number of edges}\\
\hline
$(s +2 , s +2)$ & $r$ \\
\hline
$(s+2, s+r-1)$ & $rs$ \\
\hline
$(r+s-1, r+s-1)$ & $C(s,2)$\\
\hline
\end{tabular}
\label{table:edges}
\end{table}\\\\
\begin{alignat*}{2}
ESO(C_r+K_s)=&r(s+2+s+2)\sqrt{2(s+2)^2}+rs(s+2+r+s-1)\sqrt{(s+2)^2+(r+s-1)^2}\\
&+C(s,2)(r+s-1+r+s-1)\sqrt{2(r+s-1)^2},\\
ESO(C_r+K_s)=&2\sqrt{2}r(s+2)^2+rs(r+2s+1)\sqrt{(s+2)^2+(r+s-1)^2}\\
&+2\sqrt{2}C(s,2)(r+s-1)^2.\\
EU(C_r+K_s)=&r\sqrt{3(s+2)^2}+rs\sqrt{(s+2)^2+(r+s-1)^2+(s+2)(r+s-1)}\\
&+C(s,2)\sqrt{3(r+s-1)^2},\\
EU(C_r+K_s)=&\sqrt{3}r(s+2)+rs\sqrt{(s+2)^2+(r+s-1)^2+(s+2)(r+s-1)}\\
&+\sqrt{3}C(s,2)(r+s-1).
\end{alignat*}
\end{proof}
 \end{theorem}
 \begin{theorem}
 Elliptic Sombor index of corona product of two paths are given by 
 \begin{alignat*}{2}
 ESO(P_r\odot P_s)=\begin{cases}
 34\sqrt{2}+20\sqrt{13}, r=s=2\\\\
20\sqrt{13}+2\sqrt{2}(s+1)^2+4(s+3)\sqrt{4+(s+1)^2}\\
+2(s-2)(s+4)\sqrt{9+(s+1)^2}; s>2, r=2\\\\
70+20\sqrt{13}+8\sqrt{2}r+32\sqrt{2}(r-3)+24(r-2)\sqrt{5}; r>2, s=2\\\\
18\sqrt{2}r(s-3)+2\sqrt{2}(r-3)(s+2)^2+10\sqrt{13}r+4(s+3)\sqrt{4+(s+1)^2}\\
+2(r-1)(s+4)\sqrt{4+(s+2)^2}+2(s-2)(s+4)\sqrt{9+(s+1)^2}\\
+(r-2)(s-2)(s+5)\sqrt{9+(s+2)^2}+2(2s+3)\sqrt{(s+1)^2+(s+2)^2}\\
;r,s>2
\end{cases}
 \end{alignat*}
 \begin{proof}
 Consider the graph $P_r\odot P_s$ having $|V(P_r\odot P_s)|=rs+r$ and $|E(P_r\odot P_s)|=2rs-1$. There are four types of degrees of vertices of graph $P_r\odot P_s$ namely $2, 3, s+1,s+2$. Using edge partition listed in Table 4 and formula(1),\\
 \begin{table}[h!]
 \caption{Number of edges for different cases in $P_r \odot P_s$.}
\centering
\begin{tabular}{|c|c|c|}
\hline
\textbf{Edge partitions} & \textbf{Number of edges} & \textbf{Condition} \\
\hline
$(2, 2)$ & $r$ & $s = 2$ \\
                 & $0$ & $s > 2$ \\
\hline
$(3,3)$ & $0$ & $s = 2$ \\
                 & $r(s-3)$ & $s > 2$ \\
\hline
$(s + 1, s+ 1)$ & $1$ & $r = 2$ \\
                 & $0$ & $s> 2$ \\
\hline
$(s+2,s+2)$ & $0$ & $r= 2$ \\
                 & $n-3$ & $r > 2$ \\
\hline
$(2,3)$ & $0$ & $s = 2$ \\
                 & $2r$ & $s > 2$ \\
\hline
$(2, s + 1)$ & $4$ & \text{always} \\
                 
\hline
$(2, s+2)$ & $0$ & $r = 2$ \\
                 & $2(r-2)$ & $r > 2$ \\
\hline
$(3,s+1)$ & $0$ & $s = 2$ \\
                 & $2(s - 2)$ & $s> 2$ \\
\hline
$(3,s + 2)$ & $0$ & $r= 2$ \\
                 & $(r - 2)(s-2)$ & $r> 2$ \\
\hline
$(s+ 1, s+ 2)$ & $0$ & $r = 2$\\
                 & $2$ & $ s > 2$ \\
\hline
\end{tabular}
\label{table:edges}
\end{table}
 \textbf{Case 1} r=s=2\\
 \begin{alignat*}{2}
 ESO(P_r\odot P_s)=&r(2+2)\sqrt{2^2+2^2}+(s+1+s+1)\sqrt{2(s+1)^2}\\
 &+4(s+1+2)\sqrt{(s+1)^2+2^2},\\
 \text{Subsitute} \,\, r=s=2 \\
ESO(P_r\odot P_s)=&34\sqrt{2}+20\sqrt{13}.\\
\end{alignat*}
\textbf{Case 2}  $r=2$, $s>2$\\
\begin{alignat*}{2}
ESO(P_r\odot P_s)=&(s+1+s+1)\sqrt{2(s+1)^2}+2r(2+3)\sqrt{2^2+3^2}\\
&+4(2+s+1)\sqrt{2^2+(s+1)^2}+2(s-2)(3+s+1)\sqrt{3^2+(s+1)^2}\\
\text{Substitute} \,\, r=2\\
ESO(P_r\odot P_s)=&20\sqrt{13}+2\sqrt{2}(s+1)^2+4(s+3)\sqrt{4+(s+1)^2}\\
&+2(s-2)(s+4)\sqrt{9+(s+1)^2}.
\end{alignat*}\\
\textbf{Case 3} $s=2$,$r>2$\\
\begin{alignat*}{2}
ESO(P_r\odot P_s)=&r(2+2)\sqrt{2^2+2^2}+(r-3)(s+2+s+2)\sqrt{2(s+2)^2}\\
&+4(2+s+1)\sqrt{2^2+(s+1)^2}+2(r-2)(s+2+2)\sqrt{2^2+(s+2)^2}\\
&+(r-2)(s-2)(3+s+2)\sqrt{3^2+(s+2)^2}\\
&+2(s+1+s+2)\sqrt{(s+1)^2+(s+2)^2}\\
\text{Substitute}\,\, s=2, \text{we get}\\
ESO(P_r\odot P_s)=&70+20\sqrt{13}8\sqrt{2}r+32\sqrt{2}(r-3)+24\sqrt{5}(r-2).\\
\end{alignat*}
\textbf{Case 4} $r,s>2$ \\ 
\begin{alignat*}{2}
ESO(P_r\odot P_s)=&r(s-3)(3+3)\sqrt{3^2+3^2}+(r-3)(s+2+s+2)\sqrt{2(s+2)^2}\\
&+2r(2+3)\sqrt{2^2+3^2}+4(s+1+2)\sqrt{(s+1)^2+2^2}\\
&+2(r-1)(s+2+2)\sqrt{(s+2)^2+2^2}+2(s-2)(s+1+3)\sqrt{(s+1)^2+3^2}\\
&+(r-2)(s-2)(s+2+3)\sqrt{(s+2)^2+3^2}\\
&+2(s+1+s+2)\sqrt{(s+1)^2+(s+2)^2}.
 \end{alignat*}
 \end{proof}
 \end{theorem}
 
 \begin{theorem}
 Euler Sombor index of corona product of two paths are given by
 \begin{alignat*}{2}
 EU(P_r\odot P_s)=\begin{cases}
 7\sqrt{3}+4\sqrt{19}, r=s=2\\\\
4\sqrt{19}+(s+1)\sqrt{3}+4\sqrt{4+(s+1)^2+2(s+1)}\\
+2(s-2)\sqrt{9+(s+1)^2+3(s+1)}; s>2, r=2\\\\
4\sqrt{19}+2\sqrt{3}r+8\sqrt{2}r+32\sqrt{2}(r-3)+24(r-2)\sqrt{5}; r>2, s=2\\\\
3\sqrt{3}r(s-3)+\sqrt{3}(r-3)(s+2)+2\sqrt{19}r+4\sqrt{4+(s+1)^2+2(s+1)}\\
+2(r-1)\sqrt{4+(s+2)^2+2(s+2)}+2(s-2)\sqrt{9+(s+1)^2+3(s+1)}\\
+(r-2)(s-2)\sqrt{9+(s+2)^2+3(s+2)}\\
+2\sqrt{(s+1)^2+(s+2)^2+(s+1)(s+2)}
\end{cases}
 \end{alignat*}
 \begin{proof}
Using the edge partitions in Table 4 and formula (2)
we have \\
\textbf{Case 1}\,\, r=s=2\\
\begin{alignat*}{2}
 EU(P_r\odot P_s)=&r\sqrt{2^2+2^2+2^2}+\sqrt{3(s+1)^2}+4\sqrt{(s+1)^2+2^2+2(s+1)},\\
 \text{Subsitute}\,\, r=s=2 \\
EU(P_r\odot P_s)=&7\sqrt{3}+4\sqrt{19}.\\
\textbf{Case 2} \,\,r=2,\, s>2,\\
EU(P_r\odot P_s)=&\sqrt{3(s+1)^2}+2r\sqrt{2^2+3^2+2\times 3}\\
&+4\sqrt{2^2+(s+1)^2+2(s+1)}+2(s-2)\sqrt{3^2+(s+1)^2+3(s+1)},
\end{alignat*}
\begin{alignat*}{2}
\text{Substitute}\, \,r=2\\
EU(P_r\odot P_s)=&4\sqrt{19}+\sqrt{3}(s+1)+4\sqrt{4+(s+1)^2+2(s+1)}\\
&+2(s-2)\sqrt{9+(s+1)^2+3(s+1)}.
\textbf{Case 3}\,\, s=2,r>2\\
EU(P_r\odot P_s)=&r\sqrt{2^2+2^2+2^2}+(r-3)\sqrt{3(s+2)^2}\\
&+4\sqrt{2^2+(s+1)^2+2(s+1)}+2(r-2)\sqrt{2^2+(s+2)^2+2(s+2)}\\
&+(r-2)(s-2)\sqrt{3^2+(s+2)^2+3(s+2)}\\
&+2\sqrt{(s+1)^2+(s+2)^2+(s+1)(s+2)},\\
\text{Substitute}\,\, s=2, \text{we get}\\
EU(P_r\odot P_s)=&4\sqrt{19}+2\sqrt{37}+2\sqrt{3}r+4\sqrt{3}(r-3)+4\sqrt{7}(r-2).\\
\textbf{Case 4}\,\, r,s>2 \\ 
ESO(P_r\odot P_s)=&r(s-3)\sqrt{3^2+3^2+3^2}+(r-3)\sqrt{3(s+2)^2}\\
&+2\sqrt{2^2+3^2+2\times 3}+4\sqrt{(s+1)^2+2^2+2(s+1)}\\
&+2(r-1)\sqrt{(s+2)^2+2^2+2(s+2)}+2(s-2)\\
&\sqrt{(s+1)^2+3^2+3(s+1)}+(r-2)(s-2)\\
&\sqrt{(s+2)^2+3^2+3(s+2)}+2\sqrt{(s+1)^2+(s+2)^2+(s+1)(s+2)}.
 \end{alignat*}
 \end{proof}
   \end{theorem}

\begin{theorem}
The Elliptic Sombor index and Euler Sombor index of corona product of two cycles ($C_r \odot C_s$) are given by 
\begin{alignat*}{2}
&ESO(C_r\odot C_s)=18rs\sqrt{2}+rs(s+5)\sqrt{3^2+(s+2)^2}+2\sqrt{2}r(s+2)^2.\\
&EU(C_r\odot C_s)=3\sqrt{3}+rs\sqrt{3^2+(s+2)^2+3(s+2)}+\sqrt{3}r(s+2).
\end{alignat*}
\begin{proof}
The graph $C_r\odot C_s$ with order $rs+r$ and size $2rs+r$ has two types of vertices, such as 3-degree vertex and $s+2$ degree vertex. The edges of $C_r\odot C_s$ can be partitioned into three types.
\begin{alignat*}{2}
&ESO(C_r\odot C_s)=rs(3+3)\sqrt{3^2+3^2}+rs(s+2+3)\sqrt{(s+2)^2+3^2}\\
&+r(s+2+s+2)\sqrt{2(s+2)^2},\\
&ESO(C_r\odot C_s)=18\sqrt{2}rs+rs(s+5)\sqrt{3^2+(s+2)^2}+2\sqrt{2}r(s+2)^2.\\
&EU(C_r\odot C_s)=rs\sqrt{3^2+3^2+3^2}+rs\sqrt{(s+2)^2+3^2+3(s+2)}+r\sqrt{3(s+2)^2},\\
&EU(C_r\odot C_s)=3\sqrt{3}rs+rs\sqrt{3^2+(s+2)^2+3(s+2)}+\sqrt{3}r(s+2).
\end{alignat*}
\end{proof}
\end{theorem}
\section{Conclusion}
In this article, elliptic Sombor and Euler Sombor indices of join and corona product of certain graphs are obtained for standard graphs such as path, cycle and complete graphs. These results can be further extend for corona product and join of other graphs.

\end{document}